\documentclass[a4paper,11pt]{amsart}
\usepackage{amssymb}
\usepackage{amsmath}
\usepackage{amsthm}
\usepackage[pdftex]{graphicx}
\usepackage[english]{babel}
\usepackage{array}
\usepackage{titletoc}
\usepackage[T1]{fontenc}  
\usepackage[applemac]{inputenc} 
\usepackage[all]{xy}
\usepackage[vcentermath]{youngtab}
\usepackage{stmaryrd}

 \numberwithin{equation}{section}

  \theoremstyle{plain}
\newtheorem{theoreme}[equation]{Theorem}
\newtheorem{lemme}[equation]{Lemma}
\newtheorem{lemme-def}[equation]{Lemma-Definition}
\newtheorem{proposition}[equation]{Proposition}
\newtheorem{corollaire}[equation]{Corollary}
  \theoremstyle{definition}
\newtheorem{definition}[equation]{Definition}

  \theoremstyle{remark}
\newtheorem{remarque}[equation]{Remark}
\newtheorem{notations}[equation]{Notations}
\newtheorem{exemple}[equation]{Example}

 \newcommand{\theo}{\begin{theoreme}}
 \newcommand{\defi}{\begin{definition}}
 \newcommand{\rema}{\begin{remarque}}
 \newcommand{\prop}{\begin{proposition}}
 \newcommand{\coro}{\begin{corollaire}}
 \newcommand{\lemm}{\begin{lemme}}
 \newcommand{\exem}{\begin{exemple}}
 \newcommand{\nota}{\begin{notations}}

 \newcommand{\etheo}{\end{theoreme}}
 \newcommand{\edefi}{\end{definition}}
 \newcommand{\erema}{\end{remarque}}
 \newcommand{\eprop}{\end{proposition}}
 \newcommand{\ecoro}{\end{corollaire}}
 \newcommand{\elemm}{\end{lemme}}
 \newcommand{\eexem}{\end{exemple}}
 \newcommand{\enota}{\end{notations}}

 \newcommand{\becen}{\begin{center}}
 \newcommand{\ecen}{\end{center}}

 \newcommand{\benu}{\begin{enumerate}}
 \newcommand{\eenu}{\end{enumerate}}

 \newcommand{\bite}{\begin{itemize}}
 \newcommand{\eite}{\end{itemize}}

\def\[#1\]{\begin{align*}#1\end{align*}}

\newcommand{\demo}{\begin{proof}}
\newcommand{\edemo}{\end{proof}}

\newcommand{\CC}{\mathbb C}
\newcommand{\QQ}{\mathbb Q}
\newcommand{\ZZ}{\mathbb Z}
\newcommand{\NN}{\mathbb N}

\newcommand{\Hom}{\operatorname{Hom}}
\newcommand{\End}{\operatorname{End}}

\newcommand{\aaa}{\alpha}
\newcommand{\bbb}{\beta}
\newcommand{\ccc}{\gamma}
\newcommand{\eee}{\epsilon}

\newcommand{\Irr}{\operatorname{Irr}}

\newcommand{\rarrow}{\rightarrow}

\newcommand{\Lrarrow}{\Leftrightarrow}
\newcommand{\codim}{\operatorname{codim}}

 \newcommand{\dsum}{\displaystyle\sum}
 
 \newcommand{\dprod}{\displaystyle\prod}

\newcommand{\isom}{\overset{\sim}{\rightarrow}}

\newcommand{\id}{\operatorname{id}}
\newcommand{\ima}{\operatorname{Im}}
\newcommand{\tra}{\operatorname{Tr}}
\newcommand{\W}{\textup{\sffamily W}}

\newcommand{\w}{\textup{\sffamily c}}

\author{\scshape Tristan Bozec} 

\title{ {Quivers with loops and perverse sheaves}}
\date{}

\begin{document}

\maketitle 
\tableofcontents

\section*{Introduction}
Lusztig defined in~\cite{lulu} a \emph{canonical basis} of the quantum group attached to any quiver without loop. This definition was possible thanks to an isomorphism between this quantum group and the Grothendieck group of a category of perverse sheaves, generated by the so-called \emph{Lusztig sheaves}.
Lusztig endowed this Grothendieck group with a structure of Hopf algebra, by means of restriction and induction functors.
These functors made it possible for him to perform induction proofs via a nice stratification of his category. This construction yielded a combinatorial structure on the canonical basis which would later be recognized as a \emph{Kashiwara crystal}.

There are more and more evidences of the relevance of the study of quivers with loops. A particular class of such quivers are the comet-shaped quivers, which have recently been used by Hausel, Letellier and Rodriguez-Villegas in their study of the topology of character varieties, where the number of loops at the central vertex is the genus of the considered curve (see~\cite{MR2453601} and~\cite{MR3003926}). We can also see quivers with loops appearing in a work of Nakajima relating quiver varieties with branching (see~\cite{MR2470410}), as in the work of Okounkov and Maulik about quantum cohomology (see~\cite{OM}).

Kang and Schiffmann generalized Lusztig constructions in the framework of generalized Kac-Moody algebra in~\cite{kangschiffmann}, using quivers with loops. In this case, one has to impose a somewhat unnatural restriction on the definition of a category of perverse sheaves, considering only those attached to complete flags on imaginary vertices.

In this article we consider the general definition of Lustig sheaves for arbitrary quivers, possibly carrying loops. We therefore follow the definition given in~\cite{tight}, and use the results obtained in this article for quivers with one vertex and multiple loops. Note that the category hence considered is bigger than the one considered in~\cite{kangschiffmann}, as one may already see in the case of the Jordan quiver.
We prove a conjecture raised by Lusztig in~\cite{tight}, asking if the more "simple" Lusztig perverse sheaves are enough to span the whole Grothendieck group considered. A partial proof was given in~\cite{lilin}.
Our proof  is also based on induction, still with the help of restriction and induction functors, but with non trivial first steps, consisting  in the study of quivers with one vertex but possible loops. We also need to consider regularity conditions on the support of our perverse sheaves to perform efficient restrictions at imaginary vertices.
From our proof emerges a new combinatorial structure on our generalized canonical basis, which is more general than the usual crystals, in that there are now more operators associated to a vertex with loops, as in~\cite{article1} (see~\ref{bij1}).

In a second part, we construct and study a Hopf algebra which generalizes the usual quantum groups. The geometric study previously made leads to a natural definition, which includes countably infinite sets of generators at imaginary roots, with higher order Serre relations and commutativity conditions imposed by the Jordan quiver case.
We prove that the positive part of this algebra is isomorphic to our Grothendieck group, thanks to the study of a nondegenerate Hopf pairing.

In a final section, we try to build a bridge with the Lagrangian varieties studied in~\cite{article1}, using our new Hopf algebra, as the classical case suggests (see~\cite{lulu}).

\subsubsection*{Acknowledgement}
I would like to thank Olivier Schiffmann for his constant support an availability during the preparation of this work.

\section{Quiver Varieties}

\subsection{Preliminaries}\label{Notations}

Let $Q$ be a quiver, with vertex set $I$ and oriented edge set $\Omega=\{h:s(h)\rarrow t(h)\}$. We will denote by $\Omega(i)$ the set of loops at $i$, and call $i$ \emph{imaginary} if $\omega_i=|\Omega(i)|≥1$, \emph{real} otherwise.

For every $\aaa=\sum_{i\in I}\aaa_i i\in\NN I$, we fix an $I$-graded vector space $V_\aaa$ of graded dimension $\aaa$.
For every $I$-graded vector space $X$, we set:\[
E_X&=\bigoplus_{h\in\Omega}\Hom(X_{s(h)},X_{t(h)}),\]
and $E_\aaa=E_{V_\aaa}$. We also denote by $G_\aaa$ the group $\prod_{i\in I}GL(V_{\aaa_ii})$, naturally acting on $E_\aaa$.
Take $m>0$ and two sequences $\mathbf i=(i_1,\ldots,i_m)$ and $\mathbf a=(a_1,\ldots,a_m)$ of $I$ and $\NN_{>0}$. We write $(\mathbf i,\mathbf a)\vdash\aaa$ if $\sum_{1≤k≤m}a_ki_k=\aaa$. We set:\begin{align*}
\mathcal F_{\mathbf i,\mathbf a}&=\left\{\W=(\{0\}=\W_0\subset\ldots\subset\W_m=V_\aaa)~\middle|~\forall k,\dim\dfrac{\W_{k}}{\W_{k-1}}=a_ki_k\right\}\\
\widetilde E_{\mathbf i,\mathbf a}&=\left\{(x,\W)\mid x_h(\W)\subseteq\W\right\}\subseteq E_\aaa\times\mathcal F_{\mathbf i,\mathbf a}\end{align*}
so that we get a proper morphism $\pi_{\mathbf i,\mathbf a}:\widetilde E_{\mathbf i,\mathbf a}\rarrow E_\aaa$ induced by the first projection.

Following~\cite{lusztigbook}, we will denote by $\mathcal M_G(X)$ the category of $G$-equivariant perverse sheaves on an algebraic variety $X$ equipped with an action of an algebraic connected group $G$.

Thanks to the decomposition theorem of Beilinson, Bernstein and Deligne (see \cite{BBD}), the complex ${\pi_{\mathbf i,\mathbf a}}_!\mathbf 1$ is semisimple. Denote by $\mathcal P_\aaa\subseteq\mathcal M_{G_\aaa}(E_\aaa)$ the additive category consisting of sums of $G_\aaa$-equivariant simple perverse sheaves appearing (possibly with a shift) in ${\pi_{\mathbf i,\mathbf a}}_!\mathbf 1$ for some $(\mathbf i,\mathbf a)\vdash\aaa$. Here $\mathbf 1$ stands for the constant perverse sheaf on $\widetilde E_{\mathbf i,\mathbf a}$.

Denote by $\mathcal Q_\aaa$ the category of complexes isomorphic to sums of shifts of sheaves of $\mathcal P_\aaa$. 

Let $\mathcal K_\aaa$ be the Grothendieck group of $\mathcal Q_\aaa$, seen as a $\ZZ[v^{\pm1}]$-module by setting $v^{\pm1}[\mathbf P]=[\mathbf P[\pm1]]$, $[\mathbf P]$ denoting the isoclass of a perverse sheaf $\mathbf P$. We will finally denote by $\mathcal B_\aaa$ the finite set of isoclasses of simple perverse sheaves in $\mathcal P_\aaa$, and we set $\mathcal B=\sqcup_\aaa\mathcal B_\aaa$.

For every $I$-graded subspace $W\subseteq V_\aaa$ of dimension $\bbb$ and codimension $\ccc$, equipped with two $I$-graded isomorphisms $p:W\isom V_{\bbb}$ and $q:V_\aaa/W\isom V_\ccc$, we have the following diagram:
\[\xymatrix{E_{\bbb}\times E_{\gamma}&E_\aaa(W)\ar[l]_{~~~\kappa}\ar[r]^{~~~\iota}&E_\aaa
}\] 
where $E_\aaa(W)=\{x\in E_\aaa\mid x(W)\subseteq W\}$, $\kappa:x\mapsto(p_*(x_W),q_*(x_{V_\aaa/W}))$ and $\iota$ is the inclusion. Note that $\kappa$ is a vector bundle.

We will also consider:
\[\xymatrix{E_\bbb\times E_{\gamma}&E^\dagger_{\bbb,\gamma}\ar[l]_{~~~~~p_1}\ar[r]^{p_2}&E_{\bbb,\gamma}\ar[r]^{p_3}&E_\aaa
}\] 
where: \begin{align*}
E^\dagger_{\bbb,\gamma}&=\left\{(x,W,r,\bar r)~\middle|~\begin{aligned}& x\in E_\aaa \\ &W\subseteq V_\aaa \text{ is $I$-graded and $x$-stable} \\& r:W\isom V_\bbb \\ &\bar r:V_\aaa/W\isom V_\gamma  \end{aligned}
\right\}\\
E_{\bbb,\gamma}&=\left\{(x,W)~\middle|~\begin{aligned}& x\in E_\aaa \\ &W\subseteq V_\aaa \text{ is $I$-graded and $x$-stable} \end{aligned}
\right\}.\end{align*}
These diagrams induce (\textit{cf.}~\cite[§9.2]{lusztigbook}):
\begin{align*}
\widetilde{\text {Res}}_{\bbb,\gamma}=\kappa_!\iota^*&:\mathcal Q_\aaa\rarrow\mathcal Q_\ccc\boxtimes\mathcal Q_\bbb\\
\widetilde{\textup{Ind}}_{\bbb,\ccc}={p_3}_!{p_2}_\flat p_1^*&:Q_\ccc\boxtimes\mathcal Q_\bbb\rarrow\mathcal \mathcal Q_\aaa\end{align*}
and:
\begin{align*}
\textup{Res}_{\bbb,\ccc}&=\widetilde{\text {Res}}_{\bbb,\ccc}^\aaa[d_1-d_2-2\langle \bbb, \ccc\rangle]\\
\textup{Ind}_{\bbb,\ccc}&=\widetilde{\textup{Ind}}_{\bbb,\ccc}^\aaa[d_1-d_2]\end{align*}
where $d_1$ and $d_2$ denote the dimensions of the fibers of $p_1$ and $p_2$, and $\langle\bbb,\ccc\rangle=\sum_{i\in I}\bbb_i\ccc_i$.
These functors endow $\mathcal K=\oplus_\aaa\mathcal K_\aaa$ with a Hopf algebra structure (see~\cite[10]{lulu}).
Setting $(\ccc,\bbb)=\sum_{h\in\Omega}\ccc_{s(h)}\bbb_{t(h)}$, observe that:
\begin{align*}
d_1-d_2&=(\ccc,\bbb)+\langle\bbb,\ccc\rangle\\
d_1-d_2-2\langle\bbb,\ccc\rangle&=(\ccc,\bbb)-\langle\bbb,\ccc\rangle.\end{align*}

\subsection{Study of an imaginary sink}

Let $i$ be an imaginary sink, and $(\mathbf i,\mathbf a)\vdash\aaa$. Take $\mathbf a_i=(a_{k_1},\ldots,a_{k_r})$ where $k_j<k_{j+1}$ and $\{k_j\}_{1≤j≤r}=\{k\mid i_k=i\}$. For $x\in E_\aaa$, we set $x^{(i)}=(x_h)_{h\in\Omega(i)}$ and $x^\diamond=(x_h)_{h\notin\Omega(i)}$. Then, we define:\[
\widetilde E^{(i)}_{\mathbf i,\mathbf a}&=\{(x,\W^{(i)})\mid x^{(i)}(\W^{(i)})\subseteq\W^{(i)}\}\subseteq E_\aaa\times\mathcal F_{\mathbf a_i}^{(i)}\\
E^\diamond_\aaa&=\{x\in E_\aaa\mid x^{(i)}=0\}
\]
where $\mathcal F_{\mathbf a_i}^{(i)}$ denotes the variety of flags of $V_{\aaa_ii}$ of dimension $\mathbf a_i$.
We have the following diagram:
\begin{align}\label{diagdebase}\xymatrix{
\widetilde E_{\mathbf i,\mathbf a} \ar@{}[dr]|{\square}\ar[d]_\psi \ar@/^2pc/[rr]^{\pi_{\mathbf i,\mathbf a}} \ar[r]^{\pi'_{\mathbf i,\mathbf a}}&\widetilde E^{(i)}_{\mathbf i,\mathbf a}  \ar[d]^{\text V_{\mathbf a_i}} \ar[r]^{\pi''_{\mathbf i,\mathbf a}}&E_\aaa\\
\widetilde E^\diamond_{\mathbf i,\mathbf a}\ar[r]_{\!\!\!\!\!\!\!\!\!\!\!\!\phi}&E^\diamond_\aaa\times\mathcal F_{\mathbf a_i}^{(i)}
}\end{align}
where $\widetilde E^\diamond_{\mathbf i,\mathbf a}=\{(x,\W)\in \widetilde E_{\mathbf i,\mathbf a}\mid x^{(i)}=0\}$. Note that $\psi$ and $\textup V_{\mathbf a_i}$ are vector bundles.

\subsubsection{A notion of regularity}

Put:\[
E_\aaa^{i,\textup{rss}}=\{x\in E_\aaa\mid x_h\text{ is regular semisimple if }h\in\Omega(i)\}.\]
 For any constructible subsets $X\subseteq E_\aaa$, $Y\subseteq\widetilde E_{\mathbf i,\mathbf a}$ and $Z\subseteq\widetilde E^{(i)}_{\mathbf i,\mathbf a}$, we put:\[
 X^{i,\text{rss}}&=X\cap E_\aaa^{i,\textup{rss}}\\
 Y^{i,\text{rss}}&=Y\cap\pi_{\mathbf i,\mathbf a}^{-1}(E_\aaa^{i,\textup{rss}})\\
 Z^{i,\text{rss}}&=Z\cap{\pi''_{\mathbf i,\mathbf a}}^{-1}(E_\aaa^{i,\textup{rss}}).\]
We also write $\rho_\aaa:E_\aaa^{i,\textup{rss}}\hookrightarrow E_\aaa$ for the open inclusion.

\prop\label{quasismallness} Let $\mathbf P$ be any simple element of $\mathcal P_\aaa$. Then $\mathbf P=\rho_{\aaa!*}\rho_\aaa^*\mathbf P$, \textit{i.e.}\ if $\mathbf P=\textup{IC}(Y,\mathfrak L)$ for some smooth irreducible subvariety $Y\subseteq E_\aaa$ and some local system $\mathfrak L$ on $Y$, then $Y^{i,\textup{rss}}\neq\varnothing$.\eprop

\demo By definition, $\mathbf P$ appears as a simple summand of $\pi''_{\mathbf i,\mathbf a!}\mathbf Q$ for some simple component $\mathbf Q\subseteq \pi'_{\mathbf i,\mathbf a!}\mathbf 1$. Since in~\ref{diagdebase} $\psi$ is a vector bundle and the square is cartesian, $\mathbf Q\subseteq \text V_{\mathbf a_i}^*\phi_!\mathbf 1$, and thus $\mathbf Q$ is of the form $\text{IC}(X,\mathfrak K)$ where $X=\textup V_{\mathbf a_i}^{-1}(Y)$ for an irreducible smooth subvariety  $Y\subseteq E^\diamond_\aaa\times\mathcal F_{\mathbf a_i}^{(i)}$, and $\mathfrak K=\textup V_{\mathbf a_i}^*\mathfrak L$ for an irreducible local system $\mathfrak L$ on $Y$. 

In the lemma below, we call \emph{quasismall} a map of algebraic varieties $\pi:X\rarrow Y$ satisfying the following property: there exist stratifications $X=\sqcup_{j\in J} X_j$, $Y=\sqcup_{j\in J} Y_j$ over a finite set $J$ containing an element $0$ such that:\benu
\item $X_0$ and $Y_0$ are dense;
\item $\pi_{|X_j}:X_j\rarrow Y_j$ is a locally trivial fibration of fiber $F_j$ if $j\neq0$;
\item $\pi_{|X_0}:X_0\rarrow Y_0$ is a finite morphism;
\item $2\dim F_j<\codim_YY_j$ if $j\neq0$.
\eenu
\lemm\label{quasismall} Let $S$ be a smooth irreducible subvariety of $E^\diamond_\aaa\times\mathcal F_{\mathbf a_i}^{(i)}$. Put $\widetilde S=\textup V_{\mathbf a_i}^{-1}(S)$ and $\bar S=\pi''_{\mathbf i,\mathbf a}(\widetilde S)$. Then the map $\pi''_{\mathbf i,\mathbf a|\widetilde S}:\widetilde S\rarrow\bar S$ is quasismall.\elemm
\begin{proof}[Proof of the lemma] Put $\widetilde S^0=\widetilde S^{i,\text{rss}}$, which is a nonempty open dense subset of $\widetilde S$. Moreover, the restriction of $\pi''_{\mathbf i,\mathbf a}$ to $\widetilde S^0$ is a finite morphism since a regular semisimple element $x_h$ for $h\in\Omega(i)$ stabilizes only finitely many flags of subspaces of $V_{\aaa_ii}$. Put $\widetilde T=\widetilde S\setminus\widetilde S^0$. To prove that $\pi''_{\mathbf i,\mathbf a|\widetilde S}:\widetilde S\rarrow\bar S$ is quasismall, it now suffices to check that:\[
\dim(\widetilde T\times_{E_\aaa}\widetilde T)<\dim\widetilde S.\] 

Let $z=(z_{h,k})$ be a $r\times r$-matrix of nonnegative integers such that $\sum_hz_{h,k}=a_k$, $\sum_kz_{h,k}=a_h$, and set:\[
(\widetilde S\times_{E_\aaa}\widetilde S)_z=\left\{(x,\W,\W')~\middle|~\forall h,k~\dim\dfrac{\W_h\cap\W'_k}{\W_{h-1}\cap\W'_k+\W_h\cap\W'_{k-1}}=z_{h,k}\right\}.\]
This yields a finite stratification $\widetilde S\times_{E_\aaa}\widetilde S=\sqcup_z(\widetilde S\times_{E_\aaa}\widetilde S)_z$. We use the same notations for $S\times_{E_\aaa^\diamond}S$ and $\widetilde T\times_{E_\aaa}\widetilde T$.
The fibers of $\textup V_{\mathbf a_i|\widetilde S}:\widetilde S\rarrow S$ being the same as those of $\widetilde E_{i,\mathbf a_i}\rarrow\mathcal F_{\mathbf a_i}^{(i)}$, we have for any $z$ as above:
\begin{multline}\label{memefibre}
\dim(\widetilde S\times_{E_\aaa}\widetilde S)_z-\dim(S\times_{E_\aaa^\diamond}S)_z\\=
\dim(\widetilde E_{i,\mathbf a_i}\times_{E_{\aaa_ii}}\widetilde E_{i,\mathbf a_i})_z-\dim(\mathcal F_{\mathbf a_i}^{(i)}\times \mathcal F_{\mathbf a_i}^{(i)})_z
\end{multline}
and: \begin{multline*}
\dim(\widetilde T\times_{E_\aaa}\widetilde T)_z-\dim(S\times_{E_\aaa^\diamond}S)_z\\=
\dim(\widetilde E_{i,\mathbf a_i}\times_{U_{\aaa_ii}}\widetilde E_{i,\mathbf a_i})_z-\dim(\mathcal F_{\mathbf a_i}^{(i)}\times \mathcal F_{\mathbf a_i}^{(i)})_z
\end{multline*}
where $U_{\aaa_ii}=E_{\aaa_ii}\setminus E_{\aaa_ii}^{i,\text{rss}}$.
If $\omega_i=1$, it is very well known that the map $\widetilde E_{i,\mathbf a_i}\rarrow E_{\aaa_ii}$ is quasismall, with $E_{\aaa_ii}^{i,\text{rss}}$ being the only relevant stratum. Indeed, it is true if $\mathbf a_i=(1^{\aaa_i})$, and we have the following commutative diagram:\[
\xymatrix{\widetilde E_{i,(1^{\aaa_i})}\ar[rr]^f\ar[rd]_g&&E_{\aaa_ii}\\
&\widetilde E_{i,\mathbf a_i}\ar[ru]_h&}
\]
where $g$ is projective, hence $f$ quasismall implies $h$ quasismall.
 It follows that:\begin{align}\label{stricte}
\dim(\widetilde E_{i,\mathbf a_i}\times_{U_{\aaa_ii}}\widetilde E_{i,\mathbf a_i})_z<\dim\widetilde E_{i,\mathbf a_i}.\end{align}
By~\cite{tight}, this strict inequality is also true if $\omega_i≥2$. Indeed, the large inequality is true for any $z$ if we replace $U_{\aaa_ii}$ by $E_{\aaa_ii}$, and, since $\dim U_{\aaa_ii}<\dim E_{\aaa_ii}$:\[
\dim(\widetilde E_{i,\mathbf a_i}\times_{U_{\aaa_ii}}\widetilde E_{i,\mathbf a_i})_z<\dim(\widetilde E_{i,\mathbf a_i}\times_{E_{\aaa_ii}}\widetilde E_{i,\mathbf a_i})_z≤\dim\widetilde E_{i,\mathbf a_i},\]
hence~\ref{stricte} is still satisfied.
But then:\[
&\dim\widetilde S-\dim(\widetilde T\times_{E_\aaa}\widetilde T)_z\\
&=\dim\widetilde S-\dim(S\times_{E_\aaa^\diamond}S)_z+\dim(S\times_{E_\aaa^\diamond}S)_z-\dim(\widetilde S\times_{E_\aaa}\widetilde S)_z\\
&=\dim\widetilde S-\dim(S\times_{E_\aaa^\diamond}S)_z\\
&\qquad\qquad\qquad\qquad-\dim(\widetilde E_{i,\mathbf a_i}\times_{E_{\aaa_ii}}\widetilde E_{i,\mathbf a_i})_z+\dim(\mathcal F_{\mathbf a_i}^{(i)}\times \mathcal F_{\mathbf a_i}^{(i)})_z\\
&\qquad\qquad[\text{use }\ref{memefibre}]\\
&>\dim\widetilde S-\dim(S\times_{E_\aaa^\diamond}S)_z-\dim \widetilde E_{i,\mathbf a_i}+\dim(\mathcal F_{\mathbf a_i}^{(i)}\times \mathcal F_{\mathbf a_i}^{(i)})_z\\&\qquad\qquad[\text{use }\ref{stricte}]\\
&=\dim S-\dim(S\times_{E_\aaa^\diamond}S)_z-\dim \mathcal F_{\mathbf a_i}^{(i)}+\dim(\mathcal F_{\mathbf a_i}^{(i)}\times \mathcal F_{\mathbf a_i}^{(i)})_z\\
&\qquad\qquad[\text{use }\ref{memefibre}\text{ with }z\text{ diagonal}]\\
&=\codim_{((E_\aaa^\diamond\times\mathcal F_{\mathbf a_i}^{(i)} )\times_{E_\aaa^\diamond} (E_\aaa^\diamond\times\mathcal F_{\mathbf a_i}^{(i)} ))_z}(S\times_{E_\aaa^\diamond}S)_z-\codim_{E_\aaa^\diamond\times\mathcal F_{\mathbf a_i}^{(i)}}S\\
&≥0,
\]
the last inequality being true thanks to the following diagram:
\[\xymatrix{
(S\times_{E_\aaa^\diamond}S)_z\ar@{^{(}->}[r]&\mathfrak X~\ar@{^{(}->}[r]\ar[d]\ar@{}[dr]|{\square}&E_\aaa^\diamond\times(\mathcal F_{\mathbf a_i}^{(i)}\times \mathcal F_{\mathbf a_i}^{(i)})_z\ar[d]^{\id\times\text{pr}_1}\\
&S~\ar@{^{(}->}[r]&E_\aaa^\diamond\times\mathcal F_{\mathbf a_i}^{(i)}}\]
The lemma is proved.
\end{proof}
\textit{End of proof of proposition~\ref{quasismallness}. }For any stratum $S\subseteq\overline Y$ for $\text{IC}(Y,\mathfrak L)$, the subvariety $\widetilde S=\textup V_{\mathbf a_i}^{-1}(S)$ is a stratum for $\mathbf Q$. By~\ref{quasismall}, the restriction of $\pi''_{\mathbf i,\mathbf a}$ to each of these strata is quasismall. By an argument identical to that in~\cite[1]{addendum}, it follows that $\pi''_{\mathbf i,\mathbf a!}\mathbf Q$ is a perverse sheaf, and that moreover any simple summand of $\pi''_{\mathbf i,\mathbf a!}\mathbf Q$ is an intermediate extension to $E_\aaa$ of a simple direct summand of $\pi''_{\mathbf i,\mathbf a!}(\text V_{\mathbf a_i}^*(\mathfrak L)_{|\widetilde S^0})$ for some irreducible local system $\mathfrak L$ on a stratum $S$. In particular, it is of the form $\text{IC}(R,\mathfrak J)$ where $R$ is an open subset of $\pi''_{\mathbf i,\mathbf a}(\widetilde S^0)$ for some $S$, and $\mathfrak J$ is an irreducible local system on $R$. The proposition follows from the fact that, by construction, $\pi''_{\mathbf i,\mathbf a}(\widetilde S^0)\subseteq E_\aaa^{i,\text{rss}}$.
\edemo

\subsubsection{A notion of invariance}

For any $x\in E_\aaa$, put $V_\aaa^\diamond=\oplus_{j\neq i}V_{\aaa_jj}$ and $\mathfrak I_i(x)=\CC\langle x\rangle.V_\aaa^\diamond$, \textit{i.e.}\ the smallest subspace of $V_\aaa$ stable by $x$ and containing $V_\aaa^\diamond$.
\defi Let us write $x\sim_i x'$ if the following holds:\benu
\item $x^\diamond=x'^\diamond$;
\item $\mathfrak I_i(x)\subseteq\cap_{h\in\Omega(i)}\ker(x_h-x'_h)$;
\item $\sum_{h\in\Omega(i)}\ima(x_h-x'_h)\subseteq\mathfrak I_i(x)$.
\eenu\edefi

\lemm $\sim_i$ is an equivalence relation.\elemm

\demo~

\bite
\item Reflexivity is obvious.
\item Symmetry: if $x\sim_ix'$, then $\mathfrak I(x')=\mathfrak I(x)$ since $\CC\langle x'^\diamond\rangle.V_\aaa^\diamond=\CC\langle x^\diamond\rangle.V_\aaa^\diamond\subseteq\mathfrak I_i(x)$ and since $x^{(i)}_{|\mathfrak I_i(x)}=x'^{(i)}_{|\mathfrak I_i(x)}$. This implies $x'\sim_ix$.
\item Transitivity: if $x\sim_ix'$ and $x'\sim_ix''$, we have $\mathfrak I_i(x)=\mathfrak I_i(x')=\mathfrak I_i(x'')$, $x^{(i)}_{|\mathfrak I_i(x)}=x'^{(i)}_{|\mathfrak I_i(x)}=x''^{(i)}_{|\mathfrak I_i(x)}$, and if $h\in\Omega(i)$:\[
\ima(x_h-x''_h)\subseteq\ima(x_h-x'_h)+\ima(x'_h-x''_h)\subseteq\mathfrak I_i(x).\]
Hence $x\sim_ix''$.
\eite
\edemo

Observe that equivalence classes are affine spaces. If $x\in E_\aaa$, then the equivalence class of $x$ is of dimension equal to $\omega_i\gamma(\aaa_i-\gamma)$ where $\omega_i=|\Omega(i)|$ and $\ccc i=\codim_{V_{\aaa}}\mathfrak I_i(x)$. 

There is a stratification $E_\aaa=\sqcup_{\ccc≥0}E_{\aaa,i,\ccc}$ where:\[
E_{\aaa,i,\ccc}=\{x\in E_\aaa\mid\codim_{V_{\aaa}}\mathfrak I_i(x)=\ccc i\}.\]
Note that $E_{\aaa,i,\ccc}$ is a union of $\sim_i$-equivalence classes. This can be made more precise as follows. Fix $\ccc≤\aaa_i$ and $W\subseteq V_\aaa$ an $I$-graded subspace of codimension $\ccc i$. Let $E_{\aaa,i,\ccc}(W)= E_{\aaa,i,\ccc}\cap E_\aaa(W)$ be the closed subvariety of $E_\aaa$ of elements $x\in E_\aaa$ such that $\mathfrak I_i(x)=W$. Then, if $P=\text{Stab}_{G_\aaa}(W)$,\[
E_{\aaa,i,\ccc}=G_\aaa\times_{P}E_{\aaa,i,\ccc}(W),\]
hence the inclusion $\iota_0:E_{\aaa,i,\ccc}(W)\hookrightarrow E_{\aaa,i,\ccc}$ induces an equivalence of categories of perverse sheaves:\[
\iota_0^*[-d]:\mathcal M_{G_\aaa}(E_{\aaa,i,\ccc})\rarrow\mathcal M_{P}(E_{\aaa,i,\ccc}(W))\]
where $d=\dim(G_\aaa/P)$.
Observe also that $E_{\aaa,i,\ccc}(W)$ is itself a union of $\sim_i$-equivalence classes. Here $\iota_0$ is a restriction of the inclusion $\iota$ introduced in~\ref{Notations}, with $\ccc i$ in place of $\ccc$.

Now, as in~\ref{Notations}, fix $I$-graded isomorphisms $W\simeq V_{\aaa-\ccc i}$ and $V_\aaa/W\simeq V_{\ccc i}$. We have a natural vector bundle map:\[
\kappa_0:E_{\aaa,i,\ccc}(W)\rarrow E_{\aaa-\ccc i,i,0}\times E_{\ccc i}\]
whose fibers are precisely the $\sim_i$-equivalence classes in $E_{\aaa,i,\ccc}(W)$. Again, $\kappa_0$ is a restriction of the vector bundle $\kappa$ introduced in~\ref{Notations}, with $\ccc i$ in place of $\ccc$. There is a fully faithful embedding:\[
\kappa_0^*[\omega_id]:\mathcal M_{G_{\aaa-\ccc i}\times G_{\ccc i}}(E_{\aaa-\ccc i,i,0}\times E_{\ccc i})\rarrow
\mathcal M_{P}(E_{\aaa,i,\ccc}(W)).\]
 We say that a perverse sheaf $\mathbf P\in \mathcal M_{G_\aaa}(E_{\aaa,i,\ccc})$ is \emph{$\sigma$-invariant} (at $i$) if $\iota_0^*[-d](\mathbf P)$ belongs to the essential image of $\kappa_0^*[\omega_id]$.

\defi Let $\mathcal P_{\aaa,i,≥\ccc}\subseteq\mathcal P$ be the set of perverse sheaves supported on $E_{\aaa,i,≥\ccc}$. The notation $\mathcal P_{\aaa,i,>\ccc}$ is defined likewise, and we set $\mathcal P_{\aaa,i,\ccc}=\mathcal P_{\aaa,i,≥\ccc}\setminus\mathcal P_{\aaa,i,>\ccc}$. The terms $\mathcal P_{\aaa,i,≤\ccc}$, $\mathcal P_{\aaa,i,<\ccc}$ are defined similarly.\edefi

We will need the following technical result:

\prop Let $\mathbf P$ be any simple element of $\mathcal P_{\aaa,i,\ccc}$. Let $m:E_{\aaa,i,\ccc}\hookrightarrow E_{\aaa,i,≥\ccc}$ be the open embedding. The perverse sheaf $m^*\mathbf P\in\mathcal M_{G_\aaa}(E_{\aaa,i,\ccc})$ is $\sigma$-invariant at $i$.\eprop

\demo 
The proof follows closely that of~\ref{quasismallness}, whose notations we keep. 
In particular $\mathbf P=\text{IC}(R,\mathfrak J)$ where $R$ is an open subset of $\pi''_{\mathbf i,\mathbf a}(\widetilde S^0)$ for some $G_\aaa$-invariant stratum $S\subseteq E_\aaa^\diamond\times\mathcal F^{(i)}_{\mathbf a_i}$.
Moreover $\mathbf P$ appears in some complex:\[
\mathbf R=j_{*!}\Big(\pi''_{\mathbf i,\mathbf a!}\big((\text{V}_{\mathbf a_i}^*\mathfrak L)_{|\widetilde S^0}\big)\Big)\]
where $j:\pi''_{\mathbf i,\mathbf a}(\widetilde S^0)\hookrightarrow E_\aaa$ is the inclusion and where $\mathfrak L$ is a certain $G_\aaa$-equivariant local system on $S$. It suffices to show that $\mathbf R$ is $\sigma$-equivariant.

Consider a stratification $S=\sqcup_{k}S(k)$ where:\[ 
S(k)=\{(x^\diamond,\W)\in S\mid\ima(x^\diamond)\cap V_{\aaa_ii}\subseteq\W_k\text{ but } \ima(x^\diamond)\cap V_{\aaa_ii}\not\subseteq\W_{k-1}\}.\]
Let $k$ be maximal such that $S(k)\neq\varnothing$. Then $S(k)$ is open and dense in $S$. Denote by $\widetilde S=\sqcup_l\widetilde S(l)$ the induced stratification of $\widetilde S$. Then $\widetilde S(k)$ is also open and dense in $\widetilde S$. Finally, set:\[
\widetilde S(k)^\square=\{(x,\W)\in\widetilde S(k)^{i,\text{rss}}\mid\mathfrak I_i(x)=\W_k\}.\]
It is easy to see that $\widetilde S(k)^\square$ is open and dense in $\widetilde S(k)$, hence in $\widetilde S$.

Put $\ccc=\sum_{l>k}{\mathbf a_i}_l$ so that $\ccc=\codim_{V_{\aaa_ii}}\W_k$ for any $\W\in\mathcal F_{\mathbf a_i}^{(i)}$. Let $W$ an $I$-graded subspace of $V_\aaa$ of codimension $\ccc i$ with fixed identifications $W\simeq V_{\aaa-\ccc i}$ and $V_\aaa/W\simeq V_{\ccc i}$. Consider the following diagram:\begin{align}\label{doublesquare}
\xymatrix{
	S(k)	&	\widetilde S(k)^\square\ar[l]_{\text{V}_{\mathbf a_i}}\ar[r]^{\pi''_{\mathbf i,\mathbf a}}\ar@{}[dr]|{\square}	&	E_{\aaa,i,\ccc}	\\
	S(k,W)\ar@{^(->}[u]^{\bar\iota_0}&	\widetilde S(k,W)^\square\ar@{}[dr]|{\square}\ar@{^(->}[u]^{\tilde\iota_0}\ar[d]^{\tilde\kappa_0}\ar[l]_{\text{V}_{\mathbf a_i}}\ar[r]^{\pi''_{\mathbf i,\mathbf a}}	&	E_{\aaa,i,\ccc}(W)\ar@{^(->}[u]^{\iota_0}\ar[d]^{\kappa_0}	\\
	&	\Xi\ar[ul]^{\exists\theta}\ar[r]^{\!\!\!\!\!\!\!\!\!\!\!\!\!\!\!\!\!\!\!\!\pi''}	&E_{\aaa-\ccc i,i,0}\times E_{\ccc i}
	}\end{align}
where:\bite
\item $S(k,W)=\{(x^\diamond,\W)\mid\W_k=W\}\cap S(k)\subseteq S(k)$;
\item $\widetilde S(k,W)^\square=\{(x,\W)\mid\W_k=W\}\cap\widetilde S(k)^\square\subseteq\widetilde S(k)^\square$;
\item $\bar\iota_0$, $\tilde\iota_0$ and $\tilde\kappa_0$ stand for maps induced by $\iota_0$ and $\kappa_0$;
\item $\pi''_{\mathbf i,\mathbf a}$ and $\text V_{\mathbf a_i}$ (improperly) stand for maps induced by $\pi''_{\mathbf i,\mathbf a}$ and $\text V_{\mathbf a_i}$;
\item $\Xi=\kappa(\widetilde S(k,W)^\square)\subseteq \widetilde E^{(i)}_{\mathbf i',\mathbf a'}\times \widetilde E^{(i)}_{\mathbf i'',\mathbf a''}$ where $(\mathbf i',\mathbf a')\vdash\aaa-\ccc i$ and $(\mathbf i'',\mathbf a'')\vdash\ccc i$ are naturally induced by $(\mathbf i,\mathbf a)$ and $k$. Note the existence of an inclusion $\theta$ making commutative the triangle appearing in the diragram.
\item $\pi''$ is the restriction of $\pi''_{\mathbf i',\mathbf a'}\times\pi''_{\mathbf i'',\mathbf a''}$ to $\Xi$.
\eite
Observe that the two rightmost squares are cartesian. This is obvious for the top square. For the bottom square, this follows from the fact that for $x\in E_{\aaa,i,\ccc}$,  a flag $\W\in\mathcal F_{\mathbf a_i}^{(i)}$ satisfying $\W_k=\mathfrak I_i(x)$ is $x$-stable if and only if it is $x'$-stable for any $x'\sim_ix$.

Because $\widetilde S(k)^\square$ is open and dense in $\widetilde S^0$ and $\pi''_{\mathbf i,\mathbf a|\widetilde S^0}$ is finite, we have:\[
\mathbf R=j'_{*!}\Big(\pi''_{\mathbf i,\mathbf a!}\big((\text{V}_{\mathbf a_i}^*\mathfrak L)_{|\widetilde S(k)^\square}\big)\Big)\]
where $j':\pi''_{\mathbf i,\mathbf a}(\widetilde S(k)^\square)\hookrightarrow E_\aaa$ is the inclusion.
Note that by construction $\mathbf R$ is a direct sum of objects in $\mathcal P_{\aaa,i,\ccc}$. We have:\[
m^*\mathbf R=j''_{*!}\Big(\pi''_{\mathbf i,\mathbf a!}\big((\text{V}_{\mathbf a_i}^*\mathfrak L)_{|\widetilde S(k)^\square}\big)\Big)\]
where now $j''$ and $m$ denote the inclusions defined by the following commmutative diagram:\[
\xymatrix{
\pi''_{\mathbf i,\mathbf a}(\widetilde S(k)^\square)\ar@{^(->}[rd]_{j''}\ar@{^(->}[rr]^{~~~j'}&& E_\aaa\\
&E_{\aaa,i,\ccc}\ar@{^(->}[ru]_m}\]

Furthermore, if $j''(W):\pi''_{\mathbf i,\mathbf a}(\widetilde S(k,W)^\square)\hookrightarrow E_{\aaa,i,\ccc}(W)$ denotes the inclusion induced by $j''$, \[
\iota_0^*m^*\mathbf R&=\iota_0^*j''_{*!}\pi''_{\mathbf i,\mathbf a!}\big((\text V_{\mathbf a_i}^*\mathfrak L)_{|\widetilde S(k)^\square}\big)\\
&=j''(W)_{*!}\iota_0^*\pi''_{\mathbf i,\mathbf a!}\big((\text V_{\mathbf a_i}^*\mathfrak L)_{|\widetilde S(k)^\square}\big)\\
&\qquad\qquad\qquad[\text{since }\iota_0^*\text{ is an equivalence of categories}]\\
&=j''(W)_{*!}\pi''_{\mathbf i,\mathbf a!}\big((\text V_{\mathbf a_i}^*\mathfrak L)_{|\widetilde S(k,W)^\square}\big)\\
&\qquad\qquad\qquad[\text{the highest rightmost square in~(\ref{doublesquare}) being cartesian}]\\
&=j''(W)_{*!}\pi''_{\mathbf i,\mathbf a!}\tilde\kappa_0^*\theta^*\big(\mathfrak L_{|S(k,W)}\big)\\
&\qquad\qquad\qquad[\text{the triangle being commutative in~(\ref{doublesquare})}]\\
&=j''(W)_{*!}\kappa_0^*\pi''_{!}\theta^*\big(\mathfrak L_{|S(k,W)}\big)\\
&\qquad\qquad\qquad[\text{the lowest rightmost square in~(\ref{doublesquare}) being cartesian}]\\
&=\kappa_0^*\lambda_{*!}\pi''_{!}\theta^*\big(\mathfrak L_{|S(k,W)}\big)
\]
where $\lambda:\pi''(\Xi)\hookrightarrow E_{\aaa-\ccc i,i,0}\times E_{\ccc i}$ is the inclusion (recall that $\kappa_0$ is a vector bundle). It follows that $m^*\mathbf R$ is $\sigma$-invariant as wanted. The proposition is proved.
\edemo

\subsection{A crystal type structure on $\mathcal B$}
We keep the same notations. In particular, $i$ is an imaginary sink and $W$ is an $I$-graded subspace of $V_\aaa$ of codimension $\ccc i$, with stabilizer $P\subseteq G_\aaa$. We also denote by $U$ the unipotent radical of $P$.

\prop\label{bij0} Set $d=\dim(G_\aaa/P)$.

\benu
\item[\textup{(1)}] Consider $A\in\mathcal P_{\aaa-\ccc i,i,0}\boxtimes\mathcal P_{\ccc i}$. For every $n$ we have:\[
\textup{supp}(H^n\textup{Ind}_{\aaa-\ccc i,\ccc i}A)\subseteq\overline {E_{\aaa,i,\gamma}}.\]
 If $n\neq0$, we have:\[
\textup{supp}(H^n\textup{Ind}_{\aaa-\ccc i,\ccc i}A)\cap E_{\aaa,i,\gamma}=\varnothing.
\]
Otherwise, the sum of the simple components of $H^{0}{\textup{Ind}}_{\aaa-\ccc i,\ccc i}A$ belonging to $\mathcal P_{\aaa,i,\gamma}$ is nontrivial, and we denote it by $\xi(A)$.
\item[\textup{(2)}]  Consider $B\in\mathcal P_{\aaa,i,\gamma}$. If $n\neq-2\omega_id$, we have:\[
\textup{supp}(H^n{\textup{Res}}_{\aaa-\ccc i,\ccc i}B)\cap E_{\aaa-\ccc i,i,0}\times E_{\ccc i}=\varnothing.\]
 Otherwise, the sum of the simple components of $H^{-2\omega_id}{\textup{Res}}_{\aaa-\ccc i,\ccc i}B$ belonging to $\mathcal P_{\aaa-\ccc i,i,0}\boxtimes\mathcal P_{\ccc i}$ is nontrivial, and we denote it by $\rho(B)$.
\item[\textup{(3)}] The functors $\xi$ and $\rho$ are equivalences of categories inverse to each other.
\eenu\eprop

\demo 
We will use the following diagram:\[\xymatrix{
{G_\aaa}{\times}_P E^\text{}_{\aaa,i,\gamma}(W)\ar[r]^{~~~~~~p_0}_{~~~~~~\sim}\ar[d]_{m_0}& E^\text{}_{\aaa,i,\gamma}\ar[d]_m& E^\text{}_{\aaa,i,\gamma}(W)\ar[l]_{\!\!\!\!\!\iota_0}\ar[r]^{\!\!\!\!\!\!\!\kappa_0}&E^\text{}_{\aaa-\ccc i,i,0}\times E^\text{}_{\ccc i}\ar[d]_{\mu}\\
{G_\aaa}{\times}_P E^\text{}_\aaa(W)\ar[r]^{~~~~~p=p_3}&E_{\aaa,i,≥\ccc}&E^\text{}_\aaa(W)\ar[l]_{ \iota}\ar[r]^{\!\!\!\!\!\!\!\!\kappa}&E^\text{}_{\aaa-\ccc i}\times E^\text{}_{\ccc i}
}\]

To prove (1), we denote by $\tilde A$ the perverse sheaf ${p_2}_\flat p_1^*A[(\omega_i+1)d]$. Therefore $\widetilde{\text{Ind}}_{\aaa-\ccc i,\ccc i}A=p_!\tilde A[-(\omega_i+1)d]$, and thus the support of $\widetilde{\text{Ind}}_{\aaa-\ccc i,\ccc i}A$ is included in the image of $p$, equal to $\overline {E_{\aaa,i,\gamma}}$. The following sheaf:\becen$
m^*\widetilde{\text{Ind}}_{\aaa-\ccc i,\ccc i}A=m^*p_!\tilde A[-(\omega_i+1)d]=p_{0!}m_0^*\tilde A[-(\omega_i+1)d]$\ecen
is perverse since $m_0$ is an open embedding, and since $p_0$ is an isomorphism.
The support of $H^n\widetilde{\text{Ind}}_{\aaa-\ccc i,\ccc i}A$ being included in $\overline {E_{\aaa,i,\gamma}}$ for all $n$, we get for $n\neq0$:\[
m^*H^n\widetilde{\text{Ind}}_{\aaa-\ccc i,\ccc i}A=H^nm^*\widetilde{\text{Ind}}_{\aaa-\ccc i,\ccc i}A=0\]
which proves (1) since $\widetilde{\text{Ind}}_{\aaa-\ccc i,\ccc i}A[(\omega_i+1)d]={\text{Ind}}_{\aaa-\ccc i,\ccc i}A$.

To prove (2), we use the fact that $m^*B$ is $\sigma$-equivariant, which implies that ${\kappa_0}_! \iota_0^*m^*B[-(\omega_i+1)d]$ is perverse. But:\begin{align*}
{\kappa_0}_! \iota_0^*m^*B[-(\omega_i+1)d]&=\mu^*\kappa_!\iota^*B[-(\omega_i+1)d]\\
&=\mu^*\widetilde{\text{Res}}_{\aaa-\ccc i,\ccc i}B[-(\omega_i+1)d],\end{align*} 
hence $\mu^*{\text{Res}}_{\aaa-\ccc i,\ccc i}B[-2\omega_id]$ is perverse. Since $\mu$ is an open embedding, we have, for $n\neq-2\omega_id$:\[
\mu^*H^n{\text{Res}}_{\aaa-\ccc i,\ccc i}B=H^n\mu^*{\text{Res}}_{\aaa-\ccc i,\ccc i}B=0\]
   which ends the proof of (2).
 
   We have the following diagram:
   \[\xymatrix{
E_{\aaa,i,\gamma}(W)\ar@{^(->}[d]&{G_\aaa}{\times} E_{\aaa,i,\gamma}(W)\ar[l]_{\text{pr}_{2,0}}\ar[r]^{\pi_0^P}\ar@{^(->}[d]&{G_\aaa}{\times}_P E_{\aaa,i,\gamma}(W)\ar@{^(->}[d]\\
E_{\aaa}(W)\ar[d]_{\kappa}&{G_\aaa}{\times} E_{\aaa}(W)\ar[d]_{\pi^U}\ar[l]_{\text{pr}_2}\ar[r]^{\pi^P}&{G_\aaa}{\times}_P E_{\aaa}(W)\\
E_{\aaa-\ccc i}\times E_{\ccc i}&{G_\aaa}\times_U E_\aaa(W)\ar[ru]_{p_2}\ar[l]_{p_1}
}\]
   where  $\kappa\text{pr}_2=p_1\pi^U$ by definition of $p_1$, hence $\text{pr}_2^*\kappa^*=\pi^{U*}p_1^*$, then $\pi^U_{\flat}\text{pr}_2^*\kappa^*=p_1^*$, then $p_{2\flat}\pi^{U}_{\flat}\text{pr}_2^*\kappa^*=p_{2\flat}p_1^*$ and thus:\[
\pi^P_{\flat}\text{pr}_2^*\kappa^*=p_{2\flat}p_1^*\]
   since $p_{2\flat}\pi^U_{\flat}=\pi^P_{\flat}$.

From the proof of (2) we have $\mu^*\rho(B)={\kappa_0}_! \iota_0^*m^*B[-(\omega_i+1)d]$, from which we get:\begin{align*}
 m_0^*\widetilde{\rho(B)}&= m_0^*{p_2}_\flat p_1^*{\rho(B)}[(\omega_i+1)d]\\
 &= m_0^*\pi^P_{\flat}\text{pr}_2^*\kappa^*\rho(B)[(\omega_i+1)d]\\
 &=\pi^P_{0\flat}\text{pr}_{2,0}^*\kappa_0^*\mu^*\rho(B)[(\omega_i+1)d]\\
 &=\pi^P_{0\flat}\text{pr}_{2,0}^*\kappa_0^*{\kappa_0}_! \iota_0^*m^*B\\
&=\pi^P_{0\flat}\text{pr}_{2,0}^*\iota_0^*m^*B.\end{align*}
But if we denote by $a,b:{G_\aaa}{\times} E_{\aaa,i,\gamma}\rarrow E_{\aaa,i,\gamma}$ the action of ${G_\aaa}$ on $E_{\aaa,i,\gamma}$ and the second projection, we have:
\begin{align*}
\pi^P_{0\flat}\text{pr}_{2,0}^*\iota_0^*m^*B&=\pi^P_{0\flat}(\id_{G_\aaa}\times \iota_0)^*b^*m^*B\\
&=\pi^P_{0\flat}(\id_{G_\aaa}\times \iota_0)^*a^*m^*B\\
&\qquad\qquad\qquad[\text{by ${G_\aaa}$-equivariance of }B]\\
&=\pi^P_{0\flat}\pi_0^{P*}p_0^*m^*B\\
&\qquad\qquad\qquad[\text{by definition of }p_0]\\
&=p_0^*m^*B.\end{align*}

From the proof of (1), we also have $m^*\xi(A)=p_{0!}m_0^*\tilde A$, from which we get:\[
  \mu^*\rho(\xi(A))&={\kappa_0}_! \iota_0^*m^*\xi( A)[-(\omega_i+1)d]\\
  &={\kappa_0}_! \iota_0^*p_{0!}m_0^*\tilde A[-(\omega_i+1)d]\\
  &={\kappa_0}_! \iota_0^*p_{0!}\pi^P_{0\flat}\text{pr}_{2,0}^*\kappa_0^*\mu^*A\]
but we have seen earlier that for ${G_\aaa}$-equivariant sheaves we have $\text{pr}_{2,0}^*\iota_0^*=\pi_0^{P*}p_0^*$, hence $\iota_0^*p_{0!}=\text{pr}_{2,0!}\pi_0^{P*}$, and thus:\[
  \mu^*\rho(\xi(A))&={\kappa_0}_! \kappa_0^*\mu^*A\\
  &=\mu^*A\]
but also: \begin{align*}
  m^*\xi(\rho(B))&=p_{0!}m_0^*\widetilde{\rho(B)}\\
  &=p_{0!}p_0^{*}m^*B\\
  &=m^*B.
\end{align*}
   We finally get (3).
   \edemo

\prop\label{bij1} With the same hyoptheses and notations:

\benu
\item Let $B$ be a simple object of $\mathcal P_{\aaa,i,\ccc}$. We have:\[
\textup{Res}_{\aaa-\ccc i,\ccc i}B\simeq (A\boxtimes C)\oplus(\oplus_{j\in\ZZ}L_j[j])\]
where $A$ is a simple object of $\mathcal P_{\aaa-\ccc i,i,0}$, $C$ a simple object of $\mathcal P_{\ccc i}$, and $L_j$ is the tensor product of an element of $\mathcal P_{\aaa-\ccc i,i,>0}$ and an element of $\mathcal P_{\ccc i}$ for all $j$.
\item Let $(A,C)$ be a pair of simple objects of $\mathcal P_{\aaa-\ccc i,i,0}\times\mathcal P_{\ccc i}$. We have:\becen$
\textup{Ind}_{\aaa-\ccc i,\ccc i}(A\boxtimes C)\simeq B\oplus(\oplus_{j\in \ZZ}L'_j[j])$\ecen
where $B$ is a simple object of $\mathcal P_{\aaa,i,\ccc}$ and $L'_j\in\mathcal P_{\aaa,i,>\gamma}$ for all $j$.
\item The maps $[B]\mapsto ([A],[C])$ and $([A],[C])\mapsto [B]$ induced by (1) and (2) are inverse bijections between $\mathcal B_{\aaa,i,\ccc}$ and $\mathcal B_{\aaa-\ccc i,i,0}\times\mathcal B_{\ccc i}$.
\eenu
\eprop

\demo As in~\cite[10.3.2]{lusztigbook}, the proof relies on~\ref{bij0}, using the Fourier-Deligne transform (the result~\cite[10.3.1]{lusztigbook} remains true in our setting).\edemo

We are now able to answer a question asked by Lusztig in \cite[7]{tight}. We put $\mathbf 1_{ai}={\pi_{i,a}}_!\mathbf1$:

\prop\label{psisurj} The elements $[\mathbf 1_{ai}]$ generate $\mathcal K$ ($i\in I$, $a\in\NN_{≥1}$).\eprop

\demo
We proceed by induction on $\aaa$. 
Let  $B$ be a simple object of $\mathcal P_\aaa$. Using the Fourier-Deligne transform, we may assume that there is a sink $i$ such that $B\in\mathcal P_{\aaa,i,\gamma}$ for some $\gamma>0$ (see~\cite[7.2]{lulu}). We then proceed by descending induction on $\gamma$. If $i$ is real, we can conclude as in~\cite[7.3]{lulu}. If $i$ is imaginary, the second part of~\ref{bij1} together with the one vertex quiver case enable us to conclude. Indeed, the case of the Jordan quiver is well known (see \textit{e.g.}~\cite{CriSch}), and the case of the quiver with one vertex and multiple loops is treated in~\cite{tight}.\edemo

\section{A generalized quantum group}

\subsection{Generators}

Let $(-,-)$  denote the symmetric Euler form on $\ZZ I$: $(i,j)$ is equal to the opposite of the number of edges of $\Omega$ between $i$ and $j$ for $i\neq j\in I$, and $(i,i)=2-2\omega_i$. We will denote by $I^\text{re}$ (resp. $I^\text{im}$) the set of real (resp. imaginary) vertices, and by $I^\text{iso}\subseteq I^\text{im}$ the set of \emph{isotropic} vertices: vertices $i$ such that $(i,i)=0$, \textit{i.e.}\ such that $\omega_i=1$.
We also set $I_\infty=(I^\text{re}\times\{1\})\cup( I^\text{im}\times\NN_{≥1})$, and $(\iota,j)=l(i,j)$ if $\iota=(i,l)\in I_\infty$ and $j\in I$.

\defi Let {\sffamily F} denote the  $\QQ(v)$-algebra generated by $(E_\iota)_{\iota\in I_\infty}$, naturally $\NN I$-graded by $\text{deg}(E_{i,l})=li$ for $(i,l)\in I_\infty$. We put $\text {\sffamily F}[A]=\{x\in\textup{\sffamily F}\mid |x|\in A\}$ for any $A\subseteq\NN I$, where, for convenience, we denote by $|x|$ the degree of an element $x$.\edefi

For $\aaa=\sum\aaa_ii\in\ZZ I$, we set:\benu
\item[$\triangleright$] ht$(\aaa)=\sum\aaa_i$ its height;
\item[$\triangleright$] $v_\aaa=\prod v_i^{\aaa_i}$ if $v_i=v^{(i,i)/2}$.
\eenu
We endow $\text{\sffamily F}\otimes\textup{\sffamily F}$ with the following multiplication:\[
(a\otimes b)(c\otimes d)=v^{(|b|,|c|)}(ac)\otimes(bd).\]
and equip {\sffamily F} with a comultiplication $\delta$ defined by:\begin{align*}
\delta(E_{i,l})&=\dsum_{t+t'=l}v_i^{tt'}E_{i,t}\otimes E_{i,t'}\end{align*}
where $(i,l)\in I_\infty$.

\prop For any family $(\nu_\iota)_{\iota\in I_\infty}$, we can endow \textup{\sffamily F} with a bilinear form $\langle-,-\rangle$ such that:\benu
\item[$\triangleright$]$\langle x,y \rangle=0$ if $|x|\neq|y|$;
\item[$\triangleright$]$\langle E_\iota,E_\iota\rangle=\nu_\iota$ for all $\iota\in I_\infty$;
\item[$\triangleright$]$\langle ab,c\rangle=\langle a\otimes b,\delta (c)\rangle$ for all $a,b,c\in\textup{\sffamily F}$.
\eenu
\eprop

\demo Strictly analogous to~\cite[Proposition 1.2.3]{lusztigbook} or~\cite[3]{ringel0}.\edemo

\nota Take $i\in I^\text{im}$ and $\w$ a composition (\textit{i.e.}\ a tuple of positive integers) or a partition (\textit{i.e.}\ a decreasing tuple of positive integers). We put $E_{i,\w}=\prod_j E_{i,\w_j}$, $ \nu_{i,\w}=\prod_j\nu_{i,\w_j}$, and $|\w|=\sum\w_j$.
\enota

\subsection{Relations}

\prop Consider $(\iota,j)\in I_\infty\times I^\text{re}$. The element:\begin{align}\label{Serre}
\dsum_{t+t'=-(\iota,j)+1}(-1)^tE_j^{(t)}E_{\iota}E_j^{(t')}\end{align}
belongs to the radical of $\langle-,-\rangle$.\eprop

\demo Analogous to~\cite[Proposition 1.4.3]{lusztigbook} or~\cite{ringel1}.\edemo

\rema Some higher order Serre relations are studied in~\cite[Chapter 7]{lusztigbook}, where some conditions  are given to belong to the radical. However the proofs cannot be directly adapted to our setting.\erema

The following definition is motivated by the previous proposition and our knowledge of the Jordan quiver case, which is related to the classical Hall algebra (see \textit{e.g.}~\cite{HallSch}). We know that the commutators $[E_{i,l},E_{i,k}]$ lie in the radical if $i$ is isotropic.

\defi We denote by $\tilde U^+$ the quotient of \textup{\sffamily F} by the ideal spanned by the elements~\ref{Serre} and the commutators $[E_{i,l},E_{i,k}]$ for every isotropic vertex $i$, so that $\langle-,-\rangle$ is still defined on $\tilde U^+$. We denote by $ U^+$ the quotient of $\tilde U^+$ by the radical of $\langle-,-\rangle$.
 \edefi

\defi
Let $\hat U$ be the quotient of the algebra generated by $K_i^\pm$, $E_\iota$, $F_\iota$ ($i\in I$ and $\iota\in I_\infty$) subject to the following relations:\[
K_iK_j&=K_jK_i\\
K_iK_i^-&=1\\
K_jE_\iota&=v^{(j,\iota)}E_\iota K_j\\
K_jF_\iota&=v^{-(j,\iota)}F_\iota K_j\\
\sum_{t+t'=-(\iota,j)+1}(-1)^tE_j^{(t)}E_{\iota}E_j^{(t')}&=0\qquad(j\in I^\text{re})\\
\sum_{t+t'=-(\iota,j)+1}(-1)^tF_j^{(t)}F_{\iota}F_j^{(t')}&=0\qquad(j\in I^\text{re})\\
[E_{i,l},E_{i,k}]&=0\qquad(i\in I^\text{iso})\\
[F_{i,l},F_{i,k}]&=0\qquad(i\in I^\text{iso}).\]
We extend the graduation by $|K_i|=0$ and $|F_\iota|=-|E_\iota|$, and we set $K_\aaa=\prod_iK_i^{\aaa_i}$ for every $\aaa\in\ZZ I$.

We endow $\hat U$ with a comultiplication $\Delta$ defined by:\[
\Delta(K_i)&=K_i\otimes K_i\\
\Delta(E_{i,l})&=\dsum_{t+t'=l}v_i^{tt'}E_{i,t}K_{t'i}\otimes E_{i,t'}\\
\Delta(F_{i,l})&=\dsum_{t+t'=l}v_i^{-tt'}F_{i,t}\otimes K_{-ti}F_{i,t'}.\]

We extend $\langle-,-\rangle$ to the subalgebra $\hat U^{≥0}\subseteq\hat U$ spanned by $(K_i^\pm)_{i\in I}$ and $(E_\iota)_{\iota\in I_\infty}$ by setting $\langle xK_i,yK_j\rangle=\langle x,y\rangle v^{(i,j)}$ for $x,y\in\tilde U^+$.

We use the Drinfeld double process to define $\tilde U$ as the quotient of $\hat U$ by the relations:\begin{align}\label{DD}
\dsum \langle a_{(1)},b_{(2)}\rangle \omega(b_{(1)})a_{(2)}&=\dsum \langle a_{(2)},b_{(1)}\rangle a_{(1)}\omega(b_{(2)})
\end{align}
for any $a,b\in \tilde U^{≥0}$, where $\omega$ is the unique involutive automorphism of $\hat U$ mapping $E_\iota$ to $F_\iota$ and $K_i$ to $K_{-i}$, and where we use the Sweedler notation, for example $\Delta(a)=\sum a_{(1)}\otimes a_{(2)}$.

Setting $x^-=\omega(x)$ for $x\in \tilde U$, we define $\langle-,-\rangle$ on the subalgebra $\tilde U^{-}\subseteq\tilde U$ spanned by $(F_\iota)_{\iota\in I_\infty}$ by setting $\langle x,y\rangle=\langle x^-,y^-\rangle$ for any $x,y\in\tilde  U^-$. We will denote by $U^-$ (resp. $U$) the quotient of $\tilde U^-$ (resp. $\tilde U$) by the radical of $\langle-,-\rangle$ restricted to $\tilde U^-$(resp. restricted to $\tilde U^-\times\tilde U^+$).\edefi

\prop\textup{\cite{xiao}} We can define $S,S^\textup{op}:U\rarrow U^\textup{op}$ (the \emph{antipode} and the \emph{skew antipode}) such that:\begin{align*}
&\mathbf m(S\otimes\mathbf1)\Delta=\mathbf m(\mathbf1\otimes S)\Delta=\epsilon\mathbf1\\
 &\mathbf m(S^\textup{op}\otimes\mathbf1)\Delta^\textup{op}=\mathbf m(\mathbf1\otimes S^\textup{op})\Delta^\textup{op}=\epsilon\mathbf1,\end{align*}
 where $\mathbf m$ denotes the multiplication, $\epsilon$ denotes the counit, which is equal to $1$ on $U^0$, and $0$ on $U^-\times U^+$, and $\Delta^\textup{op}$ denotes the composition of $\Delta$ and $\textup{op}:U\otimes U\rarrow U\otimes U$, $x\otimes y\mapsto y\otimes x$. We also know that $S^\textup{op}=S^{-1}$.
\eprop

\subsection{The case of the quiver with one vertex and multiple loops}

\lemm\label{modulo} We have $\langle E_{i,|\w|},E_{i,\w}\rangle=v_i^{\sum_{k<j}\w_k\w_j}\nu_{i,\w}$.\elemm

\demo By induction, using the definitions.\edemo

\prop\label{nondeg} Let $i\in I$ be a nonisotropic imaginary vertex. Assume that for every $l≥1$ we have:\begin{align}\label{hypo}
\langle E_{i,l},E_{i,l}\rangle\in1+v^{-1}\NN[\![ v^{-1}]\!].\end{align}
Then, for any compositions \textup{\sffamily c} and \textup{\sffamily c}', \becen$
\langle E_{i,\textup{\sffamily c}},E_{i,\textup{\sffamily c}'}\rangle\in\delta_{\textup{\sffamily c},\textup{\sffamily c}'}+v^{-1}\NN[\![ v^{-1}]\!]$.\ecen
\eprop

\demo For clarity, we forget the indices $i$ in this proof.
Notice that by definition of $\delta$, of the multiplication on $\textup{\sffamily F}\otimes\textup{\sffamily F}$, and since $(i,i)<0$, we already have:\becen$
\langle E_{\textup{\sffamily c}},E_{\textup{\sffamily c}'}\rangle\in\NN[\![ v^{-1}]\!]$.\ecen
Hence, we can work modulo $v^{-1}$, and then, setting $\w=(\w_1,\ldots,\w_r)$, $\w'=(\w'_1,\ldots,\w'_s)$, $\tilde\w=(\w_2,\ldots,\w_r)$ and $\tilde\w'=(\w'_2,\ldots,\w'_s)$, we get:\begin{align*}
\langle E_{\textup{\sffamily c}},E_{\textup{\sffamily c}'}\rangle&=\left\langle E_{\w_1}\otimes E_{\tilde\w},\dprod_{1≤j≤s}\delta(E_{\w'_j})\right\rangle\\
&=\left\langle E_{\w_1}\otimes E_{\tilde\w},\dprod_{1≤j≤s}(E_{\w'_j}\otimes1+1\otimes E_{\w'_j})\right\rangle\text{mod }v^{-1}\\
&=\left\{\begin{aligned}&0\text{ mod }v^{-1}&&\text{if }\w'_1\neq\w_1  \\ & \langle E_{\tilde \w},E_{\tilde\w'}\rangle\text{ mod }v^{-1}&&\text{otherwise}\end{aligned}\right.
\end{align*}
the second equality coming from the definition of $\delta$, and from $(i,i)<0$; the last equality coming from the definition of the multiplication on $\textup{\sffamily F}\otimes\textup{\sffamily F}$, from $(i,i)<0$, from~\ref{modulo}, and from the hypothesis of the proposition. We end the proof by induction.
\edemo

\coro Under the assumption~\ref{hypo}, the restriction of  $\langle-,-\rangle$ to $\textup{\sffamily F}[\NN i]$ is nondegenerate.\ecoro

\nota We denote by $\mathcal C_{i,l}$ the set of compositions $\w$ (resp. partitions) such that $|\w|=l$ if $(i,i)<0$ (resp. $(i,i)=0$).\enota

\subsection{Quasi $\mathcal R$-matrix}

\prop\label{prim} For any imaginary vertex $i$ and any $l≥1$, there exists a unique element $a_{i,l}\in \textup{\sffamily F}[li]$ such that, if we set $b_{i,l}=a_{i,l}^-$, we get:\benu
\item $\langle E_{i,l}\mid l≥1\rangle=\langle a_{i,l}\mid l≥1\rangle$ and $\langle F_{i,l}\mid l≥1\rangle=\langle b_{i,l}\mid l≥1\rangle$ as algebras;
\item $\langle a_{i,l},z\rangle=\langle b_{i,l},z^-\rangle=0$ for any $z\in\langle E_{i,k}\mid k<l\rangle$;
\item $a_{i,l}-E_{i,l}\in\langle E_{i,k}\mid k<l\rangle$ and $b_{i,l}-F_{i,l}\in\langle F_{i,k}\mid k<l\rangle$;
\item $\bar a_{i,l}=a_{i,l}$ and $\bar b_{i,l}=b_{i,l}$;
\item $\Delta(a_{i,l})=a_{i,l}\otimes1+K_{li}\otimes a_{i,l}$ and $\Delta(b_{i,l})=b_{i,l}\otimes K_{-li}+1\otimes b_{i,l}$;
\item $S(a_{i,l})=-K_{-li}a_{i,l}$ and $S(b_{i,l})=-b_{i,l}K_{li}$.
\eenu
\eprop

\demo The properties 2 and 3 enable us to define $a_{i,l}$ uniquely, and imply the other ones.
\edemo

\nota Consider $i\in I^\text{im}$ and $\w\in\mathcal C_{i,l}$. We set $\tau_{i,l}=\langle a_{i,l},a_{i,l}\rangle$,  $a_{i,\w}=\prod_j a_{i,\w_j}$, and $ \tau_{i,\w}=\prod_j\tau_{i,\w_j}$. Notice that $\{a_{i,\w}\mid \w\in\mathcal C_{i,l}\}$ is a basis of $\textup{\sffamily F}[l i]$.
\enota

\defi
We denote by $\delta_{i,\w},\delta^{i,\w}:\textup{\sffamily F}\rarrow\textup{\sffamily F}$ the linear maps defined by: \[
\delta(x)&=\dsum_{\w\in\mathcal C_{i,l}}\delta_{i,\w}(x)\otimes a_{i,\w}+\text{obd}\\
\delta(x)&=\dsum_{\w\in\mathcal C_{i,l}}a_{i,\w}\otimes\delta^{i,\w}(x)+\text{obd} \]
where "obd" stands for terms of bidegree not in $\NN I\times\NN i$ in the former equality, $\NN i\times\NN I$ in the latter one.
\edefi

\prop
The maps $\delta_{i,\w}$ and $\delta^{i,\w}$ preserve the radical of $\langle-,-\rangle$.
\eprop

\demo
First consider the case where $i$ is isotropic and $x$ is a commutator $[E_{i,l},E_{i,k}]$, then we have $\delta(x)=0$, and thus $\delta_{i,\w}(x)=\delta^{i,\w}(x)=0$. Thus, we can assume that $\langle-,-\rangle$ is nondegenerate on $\textup{\sffamily F}[\NN i]$.	
Consider $x$ in the radical of $\langle-,-\rangle$. If $|\w|=l$, we have, for all $y\in\textup{\sffamily F}$: \begin{align*}
0&=\langle x,ya_{i,\w}\rangle\\
&=\langle\delta(x),y\otimes a_{i,\w}\rangle\\
&=\dsum_{|\w'|=l}\langle\delta_{i,\w'}(x)\otimes a_{i,\w'},y\otimes a_{i,\w}\rangle\\
&=\dsum_{|\w'|=l}\langle\delta_{i,\w'}(x),y\rangle\langle a_{i,\w'}, a_{i,\w}\rangle.
\end{align*}
The result comes from the nondegeneracy of the restriction of $\langle-,-\rangle$  to $\textup{\sffamily F}[\NN i]$.\edemo

\lemm\label{delta} We have:\benu
\item $\langle a_{i,l},a_{i,\w}\rangle=\delta_{(l),\w}\tau_{i,l}$;
\item $\langle a_{i,l}y,z\rangle=\tau_{i,l}\langle y,\delta^{i,l}(z)\rangle$ for any $y,z\in\textup{\sffamily F}$;
\item $\langle ya_{i,l},z\rangle=\tau_{i,l}\langle y,\delta_{i,l}(z)\rangle$ for any $y,z\in\textup{\sffamily F}$.
\eenu\elemm

\demo The first point is a direct consequence of the definition of the $a_{i,l}$, and the rest comes from it.\edemo

\defi Let $U\hat\otimes U$ be the completion of $U\otimes U$ with respect to the following sequence ($t≥1$):\becen
$\mathcal F_t=\bigg(U^+U^0\dsum_{|\aaa|≥t}U^-[\aaa]\bigg)\otimes U+U\otimes \bigg(U^-U^0\dsum_{|\aaa|≥t}U^+[\aaa]\bigg)$.\ecen
\edefi

\prop For any $\aaa\in\NN I$, let $B_\aaa$ be a basis of $U^+[\aaa]=\{x\in U^+,|x|=\aaa\}$, and $\{b^*|b\in B_\aaa\}$ the dual basis with respect to $\langle-,-\rangle$. Set:\becen$
\Theta_\aaa=\dsum_{b\in B_\aaa}b^-\otimes b^*$.\ecen
Then, the element $\Theta=\sum\Theta_\aaa\in U\hat\otimes U$ satisfies:\becen$
\Delta(u)\Theta=\Theta\bar\Delta( u)$ for all $u\in U$\ecen
where $\bar\Delta(u)=\overline{\Delta(\overline u)}$ if $u\mapsto\overline{u}$ denotes the unique involutive $\QQ$-morphism of $U$ stabilizing $E_\iota$ and $F_\iota$, and mapping $K_i$ to $K_{-i}$, and $v$ to $v^{-1}$.
\eprop

\demo It's enough to check the relation on generators. For those of real degree, the proof is identical to the one of~\cite[Theorem 4.1.2]{lusztigbook}. Consider $i\in I^\text{im}$ and $l≥1$.
We have:
\begin{multline*}
\Delta(a_{i,l})\Theta=\Theta\bar\Delta( {a_{i,l}})\Lrarrow
\dsum_{b\in B}\{a_{i,l}b^-\otimes b^*+K_{li}b^-\otimes a_{i,l}b^*\\
-b^-a_{i,l}\otimes b^*-b^-K_{-li}\otimes b^*a_{i,l}\}=0
\end{multline*}
\begin{multline*}
\Lrarrow\forall z  \in U^+,    
\dsum_{b\in B}\{a_{i,l}b^-\langle b^*,z\rangle+K_{li}b^-\langle a_{i,l}b^*,z\rangle\\
-b^-a_{i,l}\langle b^*,z\rangle-b^-K_{-li}\langle b^*a_{i,l},z\rangle\}=0
\end{multline*}
\begin{multline*}
\Lrarrow\forall z  \in U^+ ,    
\dsum_{b\in B}\{a_{i,l}b^-\langle b^*,z\rangle+K_{li}b^-\tau_{i,l}\langle b^*,\delta^{i,l}(z)\rangle\\
-b^-a_{i,l}\langle b^*,z\rangle-b^-K_{-li}\tau_{i,l}\langle b^*,\delta_{i,l}(z)\rangle\}=0
\end{multline*}
\begin{align*}
\Lrarrow   \forall z \in U^+ ,~
a_{i,l}z^-+\tau_{i,l}K_{li}\delta^{i,l}(z)^-
=z^-a_{i,l}+\tau_{i,l}\delta_{i,l}(z)^-K_{-li}
\end{align*}
which is the relation~(\ref{DD}) with $a,b=a_{i,l},z$. The equivalence before the last one comes from~\ref{delta}.
The computations are the same for $U^{≤0}$:
\begin{multline*}
\Delta(b_{i,l})\Theta=\Theta\bar\Delta( {b_{i,l}})\Lrarrow
\dsum_{b\in B}\{b_{i,l}b^-\otimes K_{-li} b^*+b^-\otimes b_{i,l}b^*\\
-b^-b_{i,l}\otimes b^*K_{li}-b^-\otimes b^*b_{i,l}\}=0
\end{multline*}
\begin{multline*}
\Lrarrow\forall z  \in U^+ ,    
\dsum_{b\in B}\{\langle a_{i,l}b,z\rangle K_{-li} b^*+\langle b,z\rangle b_{i,l}b^*\\
-\langle ba_{i,l},z\rangle b^*K_{li}-\langle b,z\rangle b^*b_{i,l}\}=0
\end{multline*}
\begin{multline*}
\Lrarrow\forall z  \in U^+ ,    
\dsum_{b\in B}\{\tau_{i,l}\langle b,\delta^{i,l}(z)\rangle K_{-li} b^*+\langle b,z\rangle b_{i,l}b^*\\
-\tau_{i,l}\langle b,\delta_{i,l}(z)\rangle b^*K_{li}-\langle b,z\rangle b^*b_{i,l}\}=0
\end{multline*}
\begin{align*}
\Lrarrow   \forall z \in U^+,~
\tau_{i,l}K_{-li}\delta^{i,l}(z)+b_{i,l}z=\tau_{i,l}\delta_{i,l}(z)K_{li}+zb_{i,l}
\end{align*}
which matches~(\ref{DD})$^-$ with $a,b=a_{i,l},z$.
\edemo

\rema As in~\cite[4.1.2]{lusztigbook}, one can prove that $\Theta$ is the only element satisfying $\Theta_0=1\otimes1$ and $\Delta(u)\Theta=\Theta\bar\Delta( u)$ for all $u\in U$.\erema
	
\subsection{Casimir operator}

\defi We denote by $\mathcal C$ the category of $U$-modules satisfying:\benu
\item $M=\oplus_{\aaa\in\ZZ I}M_\aaa$ where $M_\aaa=\{m\in M\mid \forall i,~K_im=v^{(\aaa,i)}m\}$;
\item For any $m\in M$, there exists $p≥0$ such that $xm=0$ as soon as $x\in\textup{\sffamily F}[\aaa]$ and $\text{ht}(\aaa)≥p$.\eenu
\edefi

\prop Set $\Omega_{≤p}=\mathbf m(S\otimes\mathbf1)(\sum_{\textup{ht}(\aaa)≤p}\Theta_\aaa)$, and $M\in\mathcal C$. Then, for every $m\in M$, the value of $\Omega(m)=\Omega_{≤p}(m)$ does not depend on $p$ for $p$ large enough, and we have the following identities of operators on $M$:\begin{align*}
K_i\Omega&=\Omega K_i\\
K_{-li}a_{i,l}\Omega&=K_{li}\Omega a_{i,l}\\
b_{i,l}K_{li}\Omega K_{li}&=\Omega b_{i,l}
\end{align*}
for any $i\in I$ and $l≥1$.
\eprop

\demo The computations are strictly analogous to those in~\cite[6.1.1]{lusztigbook}, thanks to the definition of $a_{i,l}$ and $b_{i,l}$ (see~\ref{prim}).
\edemo

\defi For any $\aaa\in\ZZ I$, we define a \emph{Verma module}:\becen$
M(\aaa)=\dfrac{U}{\dsum_{\iota\in I_\infty}UE_\iota+\dsum_{i\in I}U(K_i-v^{(i,\aaa)})}\in\mathcal C$.\ecen
\edefi

\prop Under the assumption~\ref{hypo}, we have $\tilde U^-\simeq U^-$.\eprop

\demo The proof follows~\cite{kacbook}, \cite{lusztigbook} and more specifically~\cite[Proposition 2.4]{SVDB}.
The maximal degrees of the primitive elements of the kernel of the map $\tilde U^-\rarrow U^-$ are the same as those of the primitive elements of:\[
\ker\bigg(\underset{(i,l)\in I_\infty}{\dsum}\bullet~ b_{i,l}:\underset{(i,l)\in I_\infty}{\bigoplus}M(-li){\rarrow}M(0)\bigg).\]
By maximality, if $\aaa$ is such a degree, we get $(\aaa,i)≥0$ for any vertex $i$. Indeed,~\cite[§2, properties 1.,2.,3.,4.]{SVDB} are still satisfied in our case, in particular the second one, thanks to the higher order Serre relations.

Let $C$ denote the $\QQ(v)$-linear map defined on  $M=\oplus_{(i,l)\in I_\infty} M(-li)$ by:\[
Cm=v^{f(\aaa)}\Omega m\text{ if }m\in M_\aaa,\]
where $f(\aaa)=(\aaa,\aaa+2\rho)$  and $\rho$ is defined by $(i,2\rho)=(i,i)$ for every $i\in I$. Notice that:\[
f(\aaa-li)-f(\aaa)+2l(i,\aaa) =l(l-1)(i,i).\]
For any $(i,l)\in I_\infty$, since $\Omega b_{i,l}=b_{i,l}\Omega K_{2li}$, we get:\begin{align*}
Cb_{i,l}m&=v^{f(\aaa-li)}\Omega b_{i,l}m\\
&=v^{f(\aaa-li)}b_{i,l}\Omega K_{2li}m\\
&=v^{f(\aaa-li)+2l(i,\aaa)}b_{i,l}\Omega m\\
&=v^{f(\aaa-li)+2l(i,\aaa)-f(\aaa)}b_{i,l}Cm\\
&=\left\{\begin{aligned}& v^{l(l-1)(i,i)}b_{i,l}Cm &&\text{if }i\in I^\text{im}   \\ &b_{i,l}Cm &&\text{if }i\in I^\text{re}  . \end{aligned}\right.
\end{align*}
Hence, if $m$ is a primitive vector of the kernel of the map ${\oplus}_{(i,l)\in I_\infty}M(-li){\rarrow}M(0)$ with $|m|=\aaa\in-\NN I$, we have:\begin{align}\label{pos}
f(\aaa)=\dsum_{1≤k≤r}l_k(l_k-1)(i_k,i_k)
\end{align}
where $\sum_{i\in I^\text{im}}\aaa_ii=\sum_{1≤k≤r}l_ki_k$.
Since $(\aaa,i)≥0$ for any real vertex $i$, we also have:\[
(\aaa,\aaa+2\rho)&=\dsum_{i\in I}\aaa_i(i,\aaa+i)\\
&=\dsum_{i\in I^\text{re}}\aaa_i(i,\aaa)+2\dsum_{i\in I^\text{re}}\aaa_i+\sum_{i\in I^\text{im}}\aaa_i(i,\aaa+i)\\
&≤2\dsum_{i\in I^\text{re}}\aaa_i+\sum_{i\in I^\text{im}}\aaa_i(i,\aaa+i).
\]
Combining with~\ref{pos}, we get:\[
\dsum_{1≤k≤r}l_k(l_k-1)(i_k,i_k)&≤2\dsum_{i\in I^\text{re}}\aaa_i+\sum_{i\in I^\text{im}}\aaa_i(i,\aaa+i)\\
&=2\dsum_{i\in I^\text{re}}\aaa_i+\sum_{i\in I^\text{im}}\aaa_i(\aaa_i+1)(i,i)+\dsum_{\substack{i\in I^\text{im}\\j\neq i}}\aaa_i\aaa_j(i,j)\]
and thus:\[
0≤2\dsum_{i\in I^\text{re}}\aaa_i+\dsum_{\substack{i\in I^\text{im}\\j\neq i}}\aaa_i\aaa_j(i,j)
+\sum_{i\in I^\text{im}}(i,i)\bigg(\aaa_i(\aaa_i+1)-\dsum_{i_k=i}l_k(l_k-1)\bigg).\]
Since $\sum_{i_k=i}l_k=-\aaa_i$, we have:\[ 
\aaa_i(\aaa_i+1)-\sum_{i_k=i}l_k(l_k-1)=|\aaa_i|(|\aaa_i|-1)-\sum_{i_k=i}l_k(l_k-1)≥0.\]
But we also have $\aaa_i≤0$, $(i,j)≤0$ when $i\neq j$, and $(i,i)≤0$ when $i$ is imaginary, hence:\[
2\dsum_{i\in I^\text{re}}\aaa_i+\dsum_{\substack{i\in I^\text{im}\\j\neq i}}\aaa_i\aaa_j(i,j)
+\sum_{i\in I^\text{im}}(i,i)\bigg(\aaa_i(\aaa_i+1)-\dsum_{i_k=i}l_k(l_k-1)\bigg)≤0.\]
Finally every term in the sum is equal to $0$, and $-\aaa$ is a sum of pairwise othogonal imaginary vertices.
Since the restriction of $\langle-,-\rangle$ to $\tilde U^-[-\NN i]$ is nondegenerate for any imaginary vertex $i$, the proof is over.
\edemo

\theo We have an isomorphism of Hopf algebras $\Psi:U_\ZZ^+\isom\mathcal K$ defined by:\[
\left\{\begin{aligned}& E_{i,a}\mapsto\mathbf [\mathbf 1_{ai}]  &&   \text{if }i\in I^\textup{im}     \\  
&E_{i}^{(a)}\mapsto\mathbf [\mathbf 1_{ai}] &&   \text{if }i\in I^\textup{re} \end{aligned}\right.\]
 and mapping $\langle-,-\rangle$ to the geometric form $\{-,-\}$.\etheo

\demo First, $\Psi$ is defined. Indeed, we know from the Jordan quiver case that the elements $(\mathbf 1_{ai})_{a≥1}$ commute if $i$ is isotropic. Moreover the higher order Serre relations are satisfied for real vertices (see~\cite[7]{lusztigbook}), and, applying the Fourier transform on the imaginary vertices, we can assume that we are working with nilpotent representations. Hence we have $\mathbf1_{ai}={\overline\QQ_l}_{\{0_a\}}$ as if there were no loops, and the higher order Serre relations are still satisfied. 
For the same reason, we know that:\[
\{\mathbf 1_{ai},\mathbf 1_{ai}\}\in1+v^{-1}\NN[\![ v^{-1}]\!].
\] 
Hence, setting $\langle E_{i,a},E_{i,a}\rangle=\{\mathbf 1_{ai},\mathbf 1_{ai}\}$, $\langle-,-\rangle$ is nondegenerate (thanks to~\ref{nondeg}). Therefore $\Psi$ is injective, and since $\Psi$ is also surjective by~\ref{psisurj}, we get the result.
\edemo

\section{Relation with constructible functions}

We denote by $\bar h:t(h)\rarrow s(h)$ the opposite arrow of $h\in\Omega$, and $\bar Q$ the quiver $(I,H=\Omega\sqcup\bar\Omega)$, where $\bar\Omega=\{\bar h\mid h\in\Omega\}$: each arrow is replaced by a pair of arrows, one in each direction, and we set  $\eee(h)=1$ if $h\in \Omega$, $\eee(h)=-1$ if $h\in \bar\Omega$. 

For any pair of $I$-graded $\CC$-vector spaces $V=(V_i)_{i\in I}$ and $V'=(V'_i)_{i\in I}$, we set:\[
\bar E(V,V')=& \bigoplus_{h\in  H}\Hom(V_{s(h)},V'_{t(h)}).\]
For any dimension vector $\aaa=(\aaa_i)_{i\in I}$, we fix an $I$-graded $\CC$-vector space $V_\aaa$ of dimension $\aaa$, and put $\bar E_\aaa=\bar E(V_\aaa,V_\aaa)$.
The space $\bar E_\aaa=\bar E(V_\aaa,V_\aaa)$ is endowed with a symplectic form:\[
\omega_\aaa(x,x')=\dsum_{h\in H}\tra(\epsilon(h)x_hx'_{\bar h})\]
which is preserved by the natural action of $G_\aaa$ on $\bar E_\aaa$. The associated moment map $\mu_\aaa:\bar E_\aaa\rarrow\mathfrak g_\aaa= \oplus_{i\in I}\End(V_\aaa)_i$ is given by:\[
\mu_\aaa(x)=\dsum_{h\in H}\epsilon(h)x_{\bar h}x_h.\]
Here we have identified $\mathfrak g_\aaa^*$ with $\mathfrak g_\aaa$ via the trace pairing.

\defi\label{nillu}
An element $x \in \bar E_\aaa$ is said to be \emph{seminilpotent} if there exists an $I$-graded flag $\W=(\W_0=\{0\}\subset\ldots\subset\W_r=V_\aaa)$ of $V_\aaa$ such that:\[
\left.\begin{array}{ll}  
x_h(\W_\bullet)\subseteq \W_{\bullet-1}&\text{if }h\in \Omega,\\
x_h(\W_\bullet)\subseteq \W_{\bullet}&\text{if }h\in \bar\Omega.\end{array}\right.\]
We put $\Lambda(\aaa)=\{x\in \mu_\aaa^{-1}(0)\mid x\text{ seminilpotent}\}$. 
\edefi

The following is proved~\cite{article1}:

\theo The subvariety $\Lambda(\aaa)$ of $\bar E_\aaa$ is Lagrangian.\etheo

Following~\cite{semicanonical}, we denote by $\mathcal M(\aaa)$ the $\QQ$-vector space of constructible functions $\Lambda(\aaa)\rarrow\QQ$, which are constant on any $G_\aaa$-orbit. Then, we set $\mathcal M=\oplus_{\aaa≥0}\mathcal M(\aaa)$ which is a graded algebra once equipped with the product $*$ defined in~\cite[2.1]{semicanonical}.

For $Z\in \Irr\Lambda(\aaa)$ and $f\in\mathcal M(\aaa)$, we put $\rho_Z(f)=c$ if $Z\cap f^{-1}(c)$ is an open dense subset of $Z$. 

If $i\in I^\text{im}$ and $(l)$ denotes the trivial composition or partition of $l$, we denote by $1_{i,l}$ the characteristic function of the associated irreducible component $Z_{i,(l)}\in\Irr\Lambda(le_i)$ (the component of elements $x$ such that $x_h=0$ for all $h\in\Omega(i)$). If $i\notin I^\text{im}$, we just denote by $1_i$ the function mapping to $1$ the only point in $\Lambda(e_i)$.

We have $1_{i,l}\in\mathcal M(le_i)$ for $i\in I^\text{im}$ and $1_i\in\mathcal M(e_i)$ for $i\notin  I^\text{im}$. We denote by $\mathcal M_\circ\subseteq\mathcal M$ the subalgebra generated by these functions.

The following was proved in~\cite{article1}:

\prop For every $Z\in \Irr\Lambda(\aaa)$, there exists $f\in\mathcal M_\circ(\aaa)$ such that $\rho_Z(f)=1$ and $\rho_{Z'}(f)=0$ if $Z'\neq Z$.\eprop

\prop There exists a surjective morphism $\Phi:U^+_{v=1}\rarrow\mathcal M_\circ$ defined by:\[
\left\{\begin{aligned}& E_{i,a}\mapsto 1_{i,l}  &&   \text{if }i\in I^\textup{im}     \\  
&E_{i}\mapsto 1_i &&   \text{if }i\in I^\textup{re}. \end{aligned}\right.\]
\eprop

\demo 
The morphism is well defined: first, the higher order Serre relations are mapped to $0$. Indeed, they are for real vertices (see~\cite[12.11]{lulu} and~\cite[chapitre 7]{lusztigbook}), and we work with semi-nilpotent representations. Hence they are still satisfied by definition of $Z_{i,(l)}\in\Irr\Lambda(le_i)$ ($x$ such that $x_h=0$ for all $h\in\Omega(i)$).
On the other hand, the commutators $[E_{i,l},E_{i,k}]$ are also mapped to $0$ if $i$ is isotropic, thanks to the following lemma:
\lemm
Let $Q$ be the Jordan quiver. We set $I=\{\circ\}$ and $1_k=1_{\circ,k}$. We have $[1_m,1_n]=0$ for all $m,n\in\NN$.
\elemm

\demo
Consider $(x,y)\in\Lambda(n+m)$, and set $V=\CC^{n+m}$. We have:\[
1_m*1_n(x,y)=\chi\left(\left\{W\in\text{Grass}_n V
~\middle|~\begin{aligned}& W\text{ is } (x,y)\text{-stable}   \\ & x_{|W}^{|W}=0   \\ & x_{|V/W}^{|V/W}=0     \end{aligned}
\right\}\right).\]
This is equal to $0$ except if $x\in\mathcal O_\lambda$, where $\lambda=(\lambda_1≥\lambda_2)$. Then:\[
1_m*1_n(x,y)=\chi\left(\left\{\bar W\in\text{Grass}_{n-\lambda_2}  \overline{{\ker} x}\mid  \bar W~ \bar y\text{-stable} 
\right\}\right)\]
where $\bar{~}$ stands for the quotient by $\ima x$. Also:\[
1_n*1_m(x,y)=\chi\left(\left\{\bar W\in\text{Grass}_{m-\lambda_2}  \overline{{\ker} x}\mid  \bar W~ \bar y\text{-stable} 
\right\}\right).\]
Since $n-\lambda_2+m-\lambda_2=\lambda_1-\lambda_2=\dim\overline{{\ker} x}$, we get the result by duality:\[
\left.\begin{array}{rcl}   \End(\overline{{\ker} x})&\isom&\End((\overline{{\ker} x})^*) \\
   \bar y&\mapsto&[\phi\mapsto\phi\circ\bar y]   . \end{array}\right.\]
\edemo
Finally, the surjectivity comes from the definition of $\mathcal M_\circ$.
\edemo

We conjecture that $\Phi$ is an isomorphism, which should be proved by comparing the two "crystal" structures on $\mathcal K$ and $\mathcal M_\circ$ given by the following sets of bijections:\[
\mathcal B_{\aaa,i,\ccc}&\isom\mathcal B_{\aaa-\ccc i,i,0}\times\mathcal B_{\ccc i}\\
\Irr\Lambda(\aaa)_{i,\ccc}&\isom\Irr\Lambda(\aaa-\ccc i)_{i,0}\times \Irr\Lambda(\ccc i),\]
the latter being obtained in~\cite{article1}.
To that end, the notion of crystal should be generalized, and results analogous to those obtained in~\cite{kashisaito} should be proved.

\bibliographystyle{alpha}
\bibliography{biblio}

\thanks{
\noindent
\\
Facult\'e des sciences d'Orsay, B\^atiment 425, \\
Universit\'e de Paris-Sud\\
 F-91405 Orsay Cedex, France,
\\ e-mail:\;\texttt{tristan.bozec@math.u-psud.fr}}

\end{document}